\documentclass[11pt]{amsart}
\usepackage{mathrsfs}
\usepackage{amssymb}
\usepackage{marvosym}
\usepackage{color}
\usepackage{amscd}
\usepackage{amsmath, amssymb}
\usepackage{amsfonts}
\usepackage[colorlinks,linkcolor=blue,citecolor=blue, pdfstartview=FitH]{hyperref}

\numberwithin{equation}{section}

\theoremstyle{plain}
\newtheorem{thm}{Theorem}[section]
\newtheorem{theorem}[thm]{Theorem}
\newtheorem{lemma}[thm]{Lemma}
\newtheorem{corollary}[thm]{Corollary}
\newtheorem{proposition}[thm]{Proposition}

\theoremstyle{definition}

\newtheorem{remark}[thm]{Remark}

\newtheorem{definition}[thm]{Definition}

\newtheorem{example}[thm]{Example}

\newtheorem{defn-thm}[thm]{Definition-Theorem}

\newtheorem*{ack}{Acknowledgments}

\newcommand{\btheorem}{\begin{theorem}}
\newcommand{\etheorem}{\end{theorem}}
\newcommand{\bproposition}{\begin{proposition}}
\newcommand{\eproposition}{\end{proposition}}
\newcommand{\bdefinition}{\begin{definition}}
\newcommand{\edefinition}{\end{definition}}
\newcommand{\bcorollary}{\begin{corollary}}
\newcommand{\ecorollary}{\end{corollary}}
\newcommand{\bproof}{\begin{proof}}
\newcommand{\eproof}{\end{proof}}
\newcommand{\bremark}{\begin{remark}}
\newcommand{\eremark}{\end{remark}}
\newcommand{\eexample}{\end{example}}
\newcommand{\bexample}{\begin{example}}

\newcommand{\elemma}{\end{lemma}}
\newcommand{\blemma}{\begin{lemma}}

\renewcommand{\phi}{\varphi}

\newcommand{\ee}{\end{eqnarray*}}
\newcommand{\be}{\begin{eqnarray*}}

\newcommand{\beq}{\begin{equation}}
\newcommand{\eeq}{\end{equation}}
\newcommand{\bd}{\begin{enumerate}}
\newcommand{\ed}{\end{enumerate}}

\begin{document}
\title{The $\partial\overline{\partial}$-Bochner formulas for holomorphic mappings between Hermitian manifolds and their applications}
\makeatletter
\let\uppercasenonmath\@gobble
\let\MakeUppercase\relax
\let\scshape\relax
\makeatother

\author{Kai Tang}
\address{Kai Tang, College of Mathematics and Computer Science, Zhejiang Normal University, Jinhua, Zhejiang, 321004, China.}
\email{kaitang001@zjnu.edu.cn}
\maketitle

\begin{abstract}
In this paper, we derive some $\partial\overline{\partial}$-Bochner formulas for holomorphic maps between Hermitian manifolds. As applications,
we prove some Schwarz lemma type estimates, rigidity and degeneracy theorems. For instance, we show that there is no non-constant holomorphic map
from a comapct Hermitian manifold with positive (resp. non-negative) $\ell$-second Ricci curvature to a Hermitian manifold with non-positive (resp. negative)
real bisectional curvature. These theorems generalize the results \cite{Ni1,Ni2} proved recently by L. Ni on K\"{a}hler manifolds to Hermitian manifolds. We also
derive an integral inequality for holomorphic map between Hermitian manifolds.
\end{abstract}

\maketitle

\section{Introduction}
There are many generalizations of the classical Schwarz Lemma on holomorphic maps
between unit balls via the work of Ahlfors, Chen-Cheng-Look, Lu, Mok-Yau, Royden, Yau,
etc (see \cite{Koba,Lu,Royden,Yau}). Here, we in particular recall Yau's general Schwarz Lemma \cite{Yau} that a holomorphic map from a complete K\"{a}hler manifold of Ricci
curvature bounded from below to a Hermitian manifold of holomorphic bisectional curvature bounded from above by a negative constant decreases distances.
More recently, there are significant progresses on this topic, which relaxed either the curvature
assumptions or K\"{a}hlerian condition, see \cite{Ni1,Ni2,Tossa,Yang1,Yang2} and references therein for more
details. In particular, L. Ni \cite{Ni1,Ni2} proved some new estimates interpolating the Schwarz Lemmata of Royden-Yau for holomorphic mappings between K\"{a}hler manifolds.
These more flexible estimates provide additional information on (algebraic) geometric aspects of compact K\"{a}hler manifolds with nonnegative holomorphic sectional curvature, nonnegative $Ric_{\ell}$ or positive $S_{\ell}$. One wonders what the results of L. Ni can be extended
or modified to apply to the Hermitian setting.

A classical differential geometric approach in proving Schwarz type inequalities for holomorphic
map makes use of Chern-Lu formula and the maximum principle arguments. Therefore, we first generalize the $\partial\overline{\partial}$-Bochner formulas  derived by L. Ni
\cite{Ni2} on K\"{a}hler manifolds to Hermitian manifolds. Let $f: (M,g)\rightarrow (N,h)$ be a holomorphic map between Hermitian manifolds. Assume that $\dim_{\mathbb{C}}M=m \leq n=\dim_{\mathbb{C}} N$. Let $\partial f(\frac{\partial}{\partial z^{\alpha}})=\sum_{i=1}^{n}f_{\alpha}^{i}\frac{\partial}{\partial\omega^{i}}$ with respect to local coordinates
$(z^{1},\cdot\cdot\cdot,z^{m})$ and $(\omega^{1},\cdot\cdot\cdot,\omega^{n})$. The Hermitian form $f^{\ast}h=A_{\alpha\overline{\beta}}dz^{\alpha}\wedge d\overline{z}^{\beta}$
with $A_{\alpha\overline{\beta}}=f_{\alpha}^{i}\overline{f_{\beta}^{j}}h_{i\overline{j}}$ is the pull-back $h$ via $f$. For the local Hermitian metric $A=(A_{\alpha\overline{\beta}})$ and $G=(g_{\alpha\overline{\beta}})$ we denote $A_{\ell}$ and $G_{\ell}$ be the upper-left $\ell\times \ell$ blocks of them.

\btheorem \label{tk1}
Let $W_{\ell}=\frac{\det (A_{\ell})}{\det (G_{\ell})}$ be the function defined in a small neighborhood of $p$,
 $1\leq \ell\leq m=\dim_{\mathbb{C}}M$. We assume that $g_{\alpha\overline{\beta}}=\delta_{\alpha\beta}$ at $p\in M$,
$h_{i\overline{j}}=\delta_{i\overline{j}}$ at $f(p)$. Also, we assume that $df(\frac{\partial}{\partial z^{\alpha}})=\sum_{i=1}\lambda_{\alpha}\delta_{i\alpha}\frac{\partial}{\partial \omega^{i}}$ with $|\lambda_{1}|\geq\cdot\cdot\cdot\geq| \lambda_{\alpha}|\geq\cdot\cdot\cdot\geq|\lambda_{m}|$ are the singular values of $\partial f:(T_{p}^{\prime}M,g)\rightarrow (T_{f(p)}^{\prime}N,h)$.
Then at $p$, for nonzero $W_{\ell}$, we have

\begin{align}
\frac{\partial^{2}}{\partial z^{\gamma}\partial\overline{z}^{\delta}}\log W_{\ell}=&\sum_{\alpha=1}^{\ell}R^{M}(\gamma,\overline{\delta},\alpha,\overline{\alpha})-
\sum_{\alpha=1}^{\ell}R^{N}(\gamma,\overline{\delta},\alpha,\overline{\alpha})\lambda_{\gamma}\overline{\lambda_{\delta}}-\sum_{\alpha=1}^{\ell}
\sum_{t=\ell+1}^{m}g_{\alpha\overline{t},\gamma}g_{t\overline{\alpha},\overline{\delta}}\\
&+\sum_{\alpha=1}^{\ell}\sum_{i=\ell+1}^{n}\frac{1}{|\lambda_{\alpha}|^{2}}(f^{i}_{\alpha\gamma}+h_{\alpha\overline{i},\gamma}\lambda_{\alpha}\lambda_{\gamma})
(\overline{f^{i}_{\alpha\delta}}+h_{i\overline{\alpha},\overline{\delta}}\overline{\lambda_{\alpha}}\overline{\lambda_{\delta}})\nonumber
\end{align}
where $R^{M}$ and $R^{N}$ are the Chern curvature tensors of $M$ and $N$ respectively.
\etheorem

\btheorem \label{tk2}
Let $U_{\ell}=\sum_{1\leq\alpha,\beta\leq\ell}g^{\alpha\overline{\beta}}A_{\alpha\overline{\beta}}$ be the function defined in a small neighborhood of $p$,
 $1\leq \ell\leq m=\dim_{\mathbb{C}}M$. We assume that $g_{\alpha\overline{\beta}}=\delta_{\alpha\beta}$ at $p\in M$,
$h_{i\overline{j}}=\delta_{i\overline{j}}$ at $f(p)$. Also, we assume that $df(\frac{\partial}{\partial z^{\alpha}})=\sum_{i=1}\lambda_{\alpha}\delta_{i\alpha}\frac{\partial}{\partial \omega^{i}}$ with $|\lambda_{1}|\geq\cdot\cdot\cdot\geq| \lambda_{\alpha}|\geq\cdot\cdot\cdot\geq|\lambda_{m}|$ are the singular values of $\partial f:(T_{p}^{\prime}M,g)\rightarrow (T_{f(p)}^{\prime}N,h)$.
Then at $p$, for nonzero $U_{\ell}$, we have
\begin{align}
\partial\overline{\partial} U_{\ell}=&(\sum_{\delta=1}^{\ell}R^{M}(\alpha,\overline{\beta},\delta,\overline{\delta})|\lambda_{\delta}|^{2}
-\sum_{\gamma=1}^{\ell}R^{N}(\alpha,\overline{\beta},\gamma,\overline{\gamma})|\lambda_{\gamma}|^{2}\lambda_{\alpha}\overline{\lambda_{\beta}})dz^{\alpha}\wedge
d\overline{z}^{\beta}\\
&+\langle\nabla V_{\ell},\nabla V_{\ell}\rangle \nonumber
\end{align}
where $\nabla$ is the induced connection on the bundle $E=T^{\prime\ast}M\otimes f^{\ast}(T^{\prime}N)$, $V_{\ell}=\sum_{\alpha=1}^{\ell}\sum_{i=1}^{n}f^{i}_{\alpha}
dz^{\alpha}\otimes e_{i}\in \Gamma(M,E)$, $e_{i}=f^{\ast}\frac{\partial}{\partial\omega^{i}}$.

\etheorem

\bremark
When $\ell=m=\dim_{\mathbb{C}}M$, we get $W_{m}=\frac{(f^{\ast}h)^{m}}{g^{m}}$. In particular, if the domain and target
manifolds have equal dimension, namely $m=n$, then $W_{m}$ involves volume forms ralated by a holomorphic map.  Similarly
$U_{m}=tr_{g}f^{\ast}h$ is the trace of $f^{\ast}h$ respect to $g$ when $\ell=m=\dim_{\mathbb{C}}M$. Note that if $(M,g)$
and $(N,h)$ are both K\"{a}hler manifolds, we can take the normal coordinates near $p$ and $f(p)$, such that $g_{\alpha\overline{\beta}}(p)=\delta_{\alpha\beta}$,
$dg_{\alpha\overline{\beta}}(p)=0$ and $h_{i\overline{j}}(f(p))=\delta_{ij}$, $dh_{i\overline{j}}(f(p))=0$. So the calculation of the above formulas
 is much simpler (see \cite{Ni1,Ni2}).
\eremark

Before we state the applications of theorem \ref{tk1} and theorem \ref{tk2} we first recall some basic notions (also see \cite{Ni2}). 
Assume that $f: (M^{m},g)\rightarrow (N^{n},h)$ is
a holomorphic map between two Hermitian manifolds. Let $\partial f: T^{\prime}M\rightarrow T^{\prime}N$ be the tangent map. Let $\wedge^{\ell}\partial f: \wedge^{\ell}T_{x}^{\prime}M\rightarrow \wedge^{\ell}T_{f(x)}^{\prime}N$ be the associated map defined as $\wedge^{\ell}\partial f(\upsilon_{1}\wedge\cdot\cdot\cdot\wedge
\upsilon_{\ell})=\partial f(\upsilon_{1})\wedge\cdot\cdot\cdot\wedge\partial f(\upsilon_{\ell})$. Define $\|\cdot\|_{0}$ as
\begin{align}
\|\wedge^{\ell}\partial f\|_{0}(x)=\sup_{\mathbf{a}=\upsilon_{1}\wedge\cdot\cdot\cdot\wedge\upsilon_{\ell}\neq 0,\mathbf{a}\in \wedge^{\ell}T_{x}^{\prime}M}
\frac{|\wedge^{\ell}\partial f(\mathbf{a})|}{|\mathbf{a}|} \nonumber
\end{align}
We assume $g_{\alpha\overline{\beta}}=\delta_{\alpha\beta}$ at $p$,
$h_{i\overline{j}}=\delta_{i\overline{j}}$ at $f(p)$, such that $df(\frac{\partial}{\partial z^{\alpha}})=\sum_{i=1}\lambda_{\alpha}\delta_{i\alpha}\frac{\partial}{\partial \omega^{i}}$ with $|\lambda_{1}|\geq\cdot\cdot\cdot\geq| \lambda_{\alpha}|\geq\cdot\cdot\cdot\geq|\lambda_{m}|$, then $\|\wedge^{\ell}\partial f\|_{0}(p)=|\lambda_{1}\cdot\cdot\cdot\lambda_{\ell}|$. It is also easy to see that $\|\partial f\|^{2}=tr_{g}f^{\ast}h=\sum_{\alpha,\beta=1}^{m}g^{\alpha\overline{\beta}}A_{\alpha\overline{\beta}}=\sum_{\alpha=1}^{m}|\lambda_{\alpha}|^{2}$. The second goal of this
paper is to prove some estimates for holomorphic maps between Hermitian manifolds.

\btheorem \label{tk3}
Let $f: (M^{m},g)\rightarrow (N^{n},h)$ be a holomorphic map between Hermitian manifolds with $M$ being a compact. Let $m\leq n$ and $\ell\leq m$.
\begin{itemize}
\item[(a)] Assume that the scalar curvature of $(M,g)$, $S(x)\geq -K$ and the $m$-first Ricci curvature of $(N,h)$, $Ric_{m}^{(1)}(x)\leq-\kappa$,
for some $K\geq 0$, $\kappa>0$. Then
\begin{align}
\frac{(f^{\ast}h)^{m}}{g^{m}}(x)\leq (\frac{K}{m\kappa})^{m}   \nonumber
\end{align}

\item[(b)] Assume that metric $g$ on $M$ is K\"{a}hler and $1\leq \ell<m$. Assume that the $\ell$-scalar curvature of $(M,g)$, $S_{\ell}(x)\geq -K$ and the $\ell$-first Ricci curvature of $(N,h)$, $Ric_{\ell}^{(1)}(x)\leq-\kappa$,
for some $K\geq 0$, $\kappa>0$. Then
\begin{align}
\|\wedge^{\ell}\partial f\|_{0}^{2}(x)\leq (\frac{K}{\ell\kappa})^{\ell}   \nonumber
\end{align}

\item[(c)] Assume that the $\ell$-second Ricci curvature of $(M,g)$, $Ric_{\ell}^{(2)}\geq-K$, and the real bisectional curvature of $(N,h)$, $\widetilde{B}(x)\leq-\kappa$,
for some $K\geq 0$, $\kappa>0$. Then
\begin{align}
\sigma_{\ell}(x)\leq\frac{\ell K}{\kappa}. \nonumber
\end{align}
where $\sigma_{\ell}(x)=\sum_{\sigma=1}^{\ell}|\lambda_{\alpha}|^{2}(x)$.
\end{itemize}
\etheorem

We will give specific definitions of these curvatures in the next section. For $\ell=1$, $1$-first Ricci curvature and  $1$-second Ricci curvature are both holomorphic sectional curvature. If $\ell=m=\dim M$, the $m$-first Ricci curvature is (first) Chern Ricci curvature and the $m$-second Ricci curvature is second Ricci curavture. For $1\leq\ell\leq \dim M$, they are the same when the metric is K\"{a}hler.
In an attempt to generalize Wu-Yau's Theorem (\cite{WY}) to the Hermitian case, Yang and Zheng \cite{YZ} introduced the concept of {\em real bisectional curvature} for Hermitian manifolds. When the metric is K\"ahler, this curvature is the same as the holomorphic sectional curvature $H$, while when the metric is not K\"ahler, the curvature condition is slightly stronger than $H$ at least algebraically. This condition also appeared in the recent work by Lee and Streets \cite{Lee} where it is referred to as ``positive (resp. negative) curvature operator".

The following is the rigidity and degeneracy results.
\btheorem \label{tk3.1}
Let $f: (M^{m},g)\rightarrow (N^{n},h)$ be a holomorphic map between Hermitian manifolds with $M$ being a compact. Let $m\leq n$ and $\ell\leq m$.
\begin{itemize}

\item[(a)] If $S^{M}\geq0$ and manifold $(N,h)$ has $Ric_{m}^{(1)}<0$, or $S^{M}>0$ and manifold $(N,h)$ has $Ric_{m}^{(1)}\leq0$, then $f$ must be
degenerate. The same result holds if $Ric^{M}\geq0$ and $S_{m}^{N}<0$, or $Ric^{M}>0$ and $S_{m}^{N}\leq0$.
\item[(b)] Assume that metric $g$ on $M$ is K\"{a}hler and $1\leq \ell<m$. If $S_{\ell}^{M}\geq 0$ and manifold $(N,h)$ has $Ric^{(1)}_{\ell}<0$, or $S_{\ell}^{M}>0$ and manifold $(N,h)$ has $Ric^{(1)}_{\ell}\leq0$, then $rank(f)<\ell$. The same result holds if $Ric_{\ell}^{M}\geq0$ and $S_{\ell}^{N}<0$, or
 $Ric_{\ell}^{M}>0$ and $S_{\ell}^{N}\leq0$.

\item[(c)] if manifold $M$ has $Ric_{\ell}^{(2)}>0$ and manifold $N$ has $\widetilde{B}^{N}\leq0$, or $Ric_{\ell}^{(2)}\geq0$ and $\widetilde{B}^{N}<0$,
then $f$ must be constant.
\end{itemize}
\etheorem

In particular, the theorem \ref{tk3.1} recovers the following corollary.

\bcorollary \label{tk4}
There is no non-constant holomorphic map from a compact Hermitian manifold with positive (resp. non-negative) holomorphic sectional curvature to
a Hermitian manifold with non-positive (resp. negative) holomorphic sectional curvature.
\ecorollary

Note that tha above corollary \ref{tk4} is also proved independently by X. Yang in \cite{Yang1,Yang2} using a different method.

As the application of theorem \ref{tk1},  we will also give an integral inequality for non-degenerate holomorphic
maps between two Hermitian manifolds without assuming any curvature condition. Precisely, we shall prove the followings:
\btheorem \label{tk5}
Let $(M^{m},g)$ and $(N^{n},h)$ be two Hermitian manifolds and $M$ be comapct. Assume that $\dim M=m\leq n=\dim N$. Then
 there exists a smooth real function $\psi$ on $M$ such that for any non-degenerate holomorphic map $f: M\rightarrow N$ there
 holds
\begin{align}
\int_{M}S_{g}e^{(m-1)\psi}g^{m}\leq m\int_{M}e^{(m-1)\psi}f^{\ast}(Ric_{m}^{(1)}(h))\wedge g^{m-1}
\end{align}
where $S_{g}$ is the Chern scalar curvature of $g$ and $Ric_{m}^{(1)}(h)$ is the $m$-first Ricci curvature of $h$.
\etheorem

\bremark
The above theorem \ref{tk5} recovers theorem 1.2 in \cite{Zhang} which is proved by Y. Zhang when $\dim M=\dim N$. The above result can be applied to prove degeneracy theorems for holomorphic
maps without assuming any pointwise curvature signs for both the domain and target
manifolds.

\eremark

\begin{ack}
The author is grateful to Professor Fangyang Zheng for constant encouragement and support. He wishes to express
his gratitude to Professor Lei Ni for many useful discussions on \cite{Ni1}, \cite{Ni2}.
\end{ack}

\section{Preliminaries}
\subsection{Curvatures in complex geometry.}

Let $(M,g)$ be a Hermitian manifold of dimension $\dim_{\mathbb{C}}M=m$, where $\omega=\omega_{g}$ is the metric form of a Hermitian metric $g$.
If $\omega$ is closed, that is $d\omega=0$, we call $g$ a K\"{a}hler metric.
In local holomorphic chart $(z_{1},\cdot\cdot\cdot,z^{m})$, we write
\begin{align}
\omega=\sqrt{-1}\sum_{i=1,j=1}^{m}g_{i\overline{j}}dz^{i}\wedge d\overline{z}^{j} \nonumber
\end{align}
Recall the curvature tensor $R=\{R_{i\overline{j}k\overline{l}}\}$ of the Chern connection is given by
\begin{align}
R_{i\overline{j}k\overline{l}}=-\frac{\partial^{2}g_{k\overline{l}}}{\partial z^{i}\partial \overline{z}^{j}}+g^{p\overline{q}}\frac{\partial g_{k\overline{q}}}{\partial z^{i}}\frac{\partial g_{p\overline{l}}}{\partial \overline{z}^{j}}   \nonumber
\end{align}
Then the (first) Chern Ricci curvature $Ric(\omega_{g})=tr_{g}R\in\Gamma(M,\wedge^{1,1}T^{\prime\ast}M)$ has components
\begin{align}
R_{i\overline{j}}=\sum_{k=1,l=1}^{m}g^{k\overline{l}}R_{i\overline{j}k\overline{l}}=-\frac{\partial^{2}\log\det(g)}{\partial z^{i}\partial \overline{z}^{j}} \nonumber
\end{align}
The second Chern Ricci curvature $Ric^{(2)}(\omega_{g})=tr_{\omega_{g}}R\in\Gamma(M,End(T^{\prime}M))$ has components
\begin{align}
Ric_{i\overline{j}}^{(2)}=\sum_{i=1,j=1}^{m}g^{i\overline{j}}R_{i\overline{j}k\overline{l}} \nonumber
\end{align}
Note that $Ric(\omega_{g})$ and $Ric^{(2)}(\omega_{g})$ are the same when $\omega_{g}$ is K\"{a}hler metric.
The Chern scalar curvature $S_{\omega}$ is given by
\begin{align}
S_{\omega}=\sum_{i,j,k,l=1}^{m}g^{i\overline{j}}g^{k\overline{l}}R_{i\overline{j}k\overline{l}} \nonumber
\end{align}
The holomorphic bisectional curvature $B(X,Y)$ for $X,Y$ in $T_{p}^{\prime}M$ at $p\in M$ is given by
\begin{align}
B(X,Y)=\frac{H(X,\overline{X},Y,\overline{Y})}{g(X,X)g(Y,Y)}\nonumber
\end{align}
The holomorphic sectional curvature $H(X)$ is denoted by
\begin{align}
H(X)=B(X,X)=\frac{R(X,\overline{X},Y,\overline{Y})}{g(X,X)^{2}} \nonumber
\end{align}
\bdefinition Let $(M^{m},g)$ be a Hermitian manifold and $R^{g}\in\Gamma(M,\wedge^{1,1}T^{\prime\ast}M\otimes End(T^{\prime}M))$ be
the Chern curvature tensor. Assume the $\ell$-dimensional subspace $\Sigma\subset T_{p}^{\prime}M$, $p\in M$, $1\leq \ell\leq m$. For any
$\upsilon\in \Sigma$, we define $Ric^{(1)}_{\ell}(p,\Sigma)(\upsilon,\overline{\upsilon})=\sum_{i=1}^{\ell}R(\upsilon,\overline{\upsilon},E_{i},\overline{E}_{i})$
with $\{E_{i}\}$ being a unitary basis of $\Sigma$. We say $Ric_{\ell}^{(1)}(p)<0$ if $Ric_{\ell}^{(1)}(p,\Sigma)<0$ for any $\ell$-dimensional subspace $\Sigma$.
We call $Ric_{\ell}^{(1)}(p)$ $\ell$-first Ricci curvature at $p$. Similarly, For any $\upsilon\in\Sigma$, we define $Ric_{\ell}^{(2)}(p,\Sigma)(\upsilon,\overline{\upsilon})=\sum_{i=1}^{\ell}R(E_{i},\overline{E}_{i},\upsilon,\overline{\upsilon})$.
We say $Ric_{\ell}^{(2)}(p)<0$ if $Ric_{\ell}^{(2)}(p,\Sigma)<0$ for any $\ell$-dimensional subspace $\Sigma$.
We call $Ric_{\ell}^{(2)}(p)$ $\ell$-second Ricci curvature at $p$. We define $S_{\ell}(p,\Sigma)=\sum_{i,j=1}^{\ell}R(E_{i},\overline{E}_{i},E_{j},\overline{E}_{j})$ with
${E_{i}}$ being a unitary basis of $\Sigma$. We say $S_{\ell}(p)<0$ if $S_{\ell}(p,\Sigma)<0$ for any $\ell$-dimensional subspace $\Sigma$. We call $S_{\ell}(p)$
$\ell$-scalar curvature at $p$.
\edefinition
Clearly, $Ric_{\ell}^{(1)}(p)<0$ or $Ric_{\ell}^{(2)}(p)<0$ implies that $S_{\ell}(p)<0$. $Ric_{\ell}^{(1)}$ and $Ric_{\ell}^{(2)}$ are both holomorphic sectional curvature
when $\ell=1$. If $\ell=m$, they are Chern Ricci curvature and second Ricci curvature respectively. If the metric $g$ is K\"{a}hler, $Ric_{\ell}^{(1)}$ and $Ric_{\ell}^{(2)}$
will be the same. In this case of K\"{a}hler manifolds, there are many studies in \cite{Ni1,Ni2,NZ} for $Ric_{\ell}$ and $S_{\ell}$. We will do more research on the above new curvature condition
on Hermitian manifolds in the future.

Let us recall the  concept of {\em real bisectional curvature} introduced in \cite{YZ}. Let $(M^{m},g)$ be a Hermitian manifold. Denote by $R$ the curvature tensor of the Chern connection. For $p\in M$, let $e=\{e_{1},\cdot\cdot\cdot,e_{m}\}$ be a unitary tangent frame at $p$, and let $a=\{a_{1},\cdot\cdot\cdot,a_{m}\}$ be non-negative constants with $|a|^{2}=a_{1}^{2}+\cdot\cdot\cdot+a_{m}^{2}>0$. Define the $real$ $bisectional$ $curvature$ of g by
\begin{align}
\widetilde{B}_{g}(e,a)=\frac{1}{|a|^{2}}\sum_{i,j=1}^{m}R_{i\overline{i}j\overline{j}}a_{i}a_{j}.
\end{align}
We will say that a Hermitian manifold $(M^{m},g)$ has {\em positive real bisectional curvature}, denoted by $\widetilde{B}_{g}>0$, if for any $p\in M$ and any unitary frame $e$ at $p$, any nonnegative constant $a=\{a_{1},\cdot\cdot\cdot,a_{m}\}$, it holds that $\widetilde{B}_{g}(e,a)>0$.

Recall that the holomorphic sectional curvature in the direction $\upsilon$ is defined by $H(\upsilon)=R_{\upsilon\overline{\upsilon}\upsilon\overline{\upsilon}}/|\upsilon|^{4}$. If we take $e$ so that $e_{1}$ is parallel to $\upsilon$, and take $a_{1}=1$, $a_{2}=\cdot\cdot\cdot=a_{m}=0$, then $\widetilde{B}$ becomes $H(\upsilon)$. So $\widetilde{B}>0$ ($\geq0,<0,or \leq0$) implies $H>0$ ($\geq0,<0,or \leq0$).
For a more detailed discussion of this, we refer the readers to \cite{YZ}.

\subsection{Gauduchon metric.}
Let $M^{m}$ be a compact Hermitian manifold. A Hermitian metric $\omega$ is called Gauduchon if
\begin{align}
\partial\overline{\partial}(\omega^{m-1})=0 \nonumber
\end{align}
For a Gauduchon metric $\omega$ and a smooth function $u$ on $M$, we easily get
\begin{align}
\int_{M}(\Delta_{\omega}u)\omega^{m}=0 \nonumber
\end{align}
where $\Delta_{\omega}u$ is the complex Laplacian defined by $\Delta_{\omega}u=tr_{\omega}(\sqrt{-1}\partial\overline{\partial}u)$.
A classical result of Gauduchon \cite{Gau} states that, for any Hermitian metric $\omega$, there is a $\psi\in C^{\infty}(M,\mathbb{R})$ (unique up to
scaling) such that $e^{\psi}\omega$ is Gauduchon.

\subsection{Non-degenerate holomorphic maps.}
Let $f: M^{m}\rightarrow N^{n}$ be a holomorphic map between two Hermitian manifolds $(m\leq n)$. If $\dim(f(M))=m$, then we say
$f$ is non-degenerate; If $\dim(f(M))<m$, we say $f$ is degenerate.

\section{$\partial\overline{\partial}$-Bochner formulas for holomorphic mappings}
In this section, we will give the proof of theorem \ref{tk1} and \ref{tk2}. Compared to the K\"{a}hler case, its calculation is more complicated in the Hermitian case.
\bproof [Proof of Theorem \ref{tk1}]
As stated in the theorem, we assume that $g_{\alpha\overline{\beta}}=\delta_{\alpha\overline{\beta}}$ at $p\in M$, $h_{i\overline{j}}=\delta_{i\overline{j}}$ at $f(p)$.
After a constant unitary change of coordinates $z$ and $\omega$, at $p$, we have $df(\frac{\partial}{\partial z^{\alpha}})=\sum_{i=1}^{n}\lambda_{\alpha}\delta_{i\alpha}\frac{\partial}{\partial \omega^{i}}$ with
$|\lambda_{1}|\geq\cdot\cdot\cdot\geq|\lambda_{\alpha}|\geq\cdot\cdot\cdot|\lambda_{m}|$. The Hermitian form
$f^{\ast}h=A_{\alpha\overline{\beta}}dz^{\alpha}\wedge d\overline{z}^{\beta}$ with $A_{\alpha\overline{\beta}}=\sum_{i,j=1}^{n}f^{i}_{\alpha}\overline{f^{j}_{\beta}}h_{i\overline{j}}$. For the local Hermitian metric $A=(A_{\alpha\overline{\beta}})$
and $G=(g_{\alpha\overline{\beta}})$ we denote $A_{\ell}$ and $G_{\ell}$ be the upper-left $\ell\times\ell$ blocks of them. To simplify notations we write
$\frac{\partial f^{i}}{\partial z^{\alpha}}$ as $f^{i}_{\alpha}$, and $\frac{\partial^{2}f^{i}}{\partial z^{\alpha}\partial z^{\gamma}}$ as $f^{i}_{\alpha\gamma}$.
We can perform the computation at $p$ and $f(p)$, where
\begin{align}
R^{M}_{\gamma\overline{\delta}\alpha\overline{\beta}}=
-g_{\alpha\overline{\beta},\overline{\delta}\gamma}+\sum_{t=1}^{m}g_{\alpha\overline{t},\gamma}g_{t\overline{\beta},\overline{\delta}}\,\,\,,
R^{N}_{i\overline{j}k\overline{l}}=
-g_{k\overline{l},\overline{j}i}+\sum_{s=1}^{n}g_{k\overline{s},i}g_{s\overline{l},\overline{j}} \nonumber
\end{align}
Hence,
\begin{align}
\frac{\partial^{2}}{\partial z^{\gamma}\partial \overline{z}^{\delta}}\log W_{\ell}=\frac{\partial^{2}}{\partial z^{\gamma}\partial\overline{z}^{\delta}}
\log\frac{\det(A_{\ell})}{\det(G_{\ell})}=\frac{\partial^{2}}{\partial z^{\gamma}\partial\overline{z}^{\delta}}\log\det(A_{\ell})-\frac{\partial^{2}}{\partial z^{\gamma}\partial\overline{z}^{\delta}}\log\det(G_{\ell})\nonumber
\end{align}
Direct calculation shows that
\begin{align}
\frac{\partial^{2}}{\partial z^{\gamma}\partial\overline{z}^{\delta}}\log\det(G_{\ell})&=\frac{\partial}{\partial z^{\gamma}}[\sum_{\alpha,\beta=1}^{\ell}(G_{\ell})^{\alpha\overline{\beta}}g_{\alpha\overline{\beta},\overline{\delta}}]\nonumber\\
&=-\sum_{\alpha,\beta,t,\lambda=1}^{\ell}(G_{\ell})^{\alpha\overline{t}}g_{\lambda\overline{t},\gamma}(G_{\ell})^{\lambda\overline{\beta}}
g_{\alpha\overline{\beta},\overline{\delta}}+\sum_{\alpha,\beta=1}^{\ell}(G_{\ell})^{\alpha\overline{\beta}}g_{\alpha\overline{\beta},\overline{\delta}\gamma}\nonumber\\
&=\sum_{\alpha=1}^{\ell}(-\sum_{t=1}^{\ell}g_{\alpha\overline{t},\gamma}g_{t\overline{\alpha},\overline{\delta}}+g_{\alpha\overline{\alpha},\overline{\delta}\gamma})\nonumber\\
&=\sum_{\alpha=1}^{\ell}(-R^{M}_{\gamma\overline{\delta}\alpha\overline{\alpha}}+\sum_{t=\ell+1}^{m}g_{\alpha\overline{t},\gamma}g_{t\overline{\alpha},\overline{\delta}})\nonumber
\end{align}
The last two lines only holds at point $p$.
\begin{align}
(\log\det(A_{\ell}))_{\overline{\delta}}&=\sum_{\alpha,\beta=1}^{\ell}\sum_{i,j=1}^{n}(A_{\ell})^{\alpha\overline{\beta}}(f^{i}_{\alpha}h_{i\overline{j}}
\overline{f^{j}_{\beta}})_{\overline{\delta}} \nonumber\\
&=\sum_{\alpha,\beta=1}^{\ell}\sum_{i,j=1}^{n}(A_{\ell})^{\alpha\overline{\beta}}
[\sum_{t=1}^{n}h_{i\overline{j},\overline{t}}f^{i}_{\alpha}\overline{f^{j}_{\beta}}\overline{f^{t}_{\delta}}+f^{i}_{\alpha}h_{i\overline{j}}\overline{f^{j}_{\beta\delta}}]\nonumber\\
&=\sum_{\alpha=1}^{\ell}\frac{1}{|\lambda_{\alpha}|^{2}}
[h_{\alpha\overline{\alpha},\overline{\delta}}|\lambda_{\alpha}|^{2}\overline{\lambda_{\delta}}+\lambda_{\alpha}\overline{f^{\alpha}_{\alpha\delta}}]\nonumber
\end{align}
Similarly,
\begin{align}
(\log\det(A_{\ell}))_{\gamma}&=\sum_{\alpha,\beta=1}^{\ell}\sum_{i,j=1}^{n}(A_{\ell})^{\alpha\overline{\beta}}[f^{i}_{\alpha\gamma}h_{i\overline{j}}\overline{f^{j}_{\beta}}+
\sum_{k=1}^{n}f^{i}_{\alpha}h_{i\overline{j},k}f^{k}_{\gamma}\overline{f^{j}_{\beta}}]\nonumber\\
&=\sum_{\alpha=1}^{\ell}\frac{1}{|\lambda_{\alpha}|^{2}}[f^{\alpha}_{\alpha\gamma}\overline{\lambda_{\alpha}}
+h_{\alpha\overline{\alpha},\gamma}|\lambda_{\alpha}|^{2}\lambda_{\gamma}]\nonumber
\end{align}
Taking second derivative, at the end restricting to $p$, we have
\begin{align}
(\log\det(A_{\ell}))_{\gamma\overline{\delta}}&=\frac{\partial}{\partial z^{\gamma}}\{\sum_{\alpha,\beta=1}^{\ell}\sum_{i,j=1}^{n}(A_{\ell})^{\alpha\overline{\beta}}
[\sum_{t=1}^{n}h_{i\overline{j},\overline{t}}f^{i}_{\alpha}\overline{f^{j}_{\beta}}\overline{f^{t}_{\delta}}+f^{i}_{\alpha}h_{i\overline{j}}\overline{f^{j}_{\beta\delta}}]\}\nonumber\\
&=\sum_{\alpha,\beta=1}^{\ell}\sum_{i,j=1}^{n}\{-\sum_{k,s=1}^{\ell}(A_{\ell})^{\alpha\overline{s}}(A_{\ell})_{k\overline{s},\gamma}(A_{\ell})^{k\overline{\beta}}
[\sum_{t=1}^{n}h_{i\overline{j},\overline{t}}f^{i}_{\alpha}\overline{f^{j}_{\beta}}\overline{f^{t}_{\delta}}+f^{i}_{\alpha}h_{i\overline{j}}\overline{f^{j}_{\beta\delta}}]\}\nonumber\\
&\,\,\,\,\,\,+\sum_{\alpha,\beta=1}^{\ell}\sum_{i,j=1}^{n}(A_{\ell})^{\alpha\overline{\beta}}[\sum_{t,p=1}^{n}h_{i\overline{j},\overline{t}p}f_{\gamma}^{p}f^{i}_{\alpha}
\overline{f^{j}_{\beta}}\overline{f^{t}_{\delta}}+\sum_{t=1}^{n}h_{i\overline{j},\overline{t}}f^{i}_{\alpha\gamma}\overline{f^{j}_{\beta}}\overline{f^{t}_{\delta}} \nonumber\\
&\,\,\,\,\,\,+f^{i}_{\alpha\gamma}h_{i\overline{j}}\overline{f^{j}_{\beta\delta}}+\sum_{q=1}^{n}f^{i}_{\alpha}h_{i\overline{j},q}f^{q}_{\gamma}\overline{f^{j}_{\beta\delta}}]\nonumber\\
&=-\sum_{\alpha,\beta=1}^{\ell}\frac{1}{|\lambda_{\alpha}|^{2}|\lambda_{\beta}|^{2}}(f^{\alpha}_{\beta\gamma}\overline{\lambda_{\alpha}}+
\lambda_{\beta}\lambda_{\gamma}\overline{\lambda_{\alpha}}h_{\beta\overline{\alpha},\gamma})(h_{\alpha\overline{\beta},\overline{\delta}}
\lambda_{\alpha}\overline{\lambda_{\beta}}\overline{\lambda_{\delta}}+\lambda_{\alpha}\overline{f^{\alpha}_{\beta\delta}})\nonumber\\
&\,\,\,\,\,\,+\sum_{\alpha=1}^{\ell}\sum_{i=1}^{n}\frac{1}{|\lambda_{\alpha}|^{2}}(h_{\alpha\overline{\alpha},\overline{\delta}\gamma}\lambda_{\gamma}|\lambda_{\alpha}|^{2}
\overline{\lambda_{\delta}}+h_{i\overline{\alpha},\overline{\delta}}f^{i}_{\alpha\gamma}\overline{\lambda_{\alpha}}\overline{\lambda_{\delta}}\nonumber\\
&\,\,\,\,\,\,+f^{i}_{\alpha\gamma}\overline{f^{i}_{\alpha\delta}}+\lambda_{\alpha}h_{\alpha\overline{j},\gamma}\lambda_{\gamma}\overline{f^{^{i}}_{\alpha\delta}}) \nonumber\\
&=\sum_{\alpha=1}^{\ell}\sum_{i=\ell+1}^{n}\frac{1}{|\lambda_{\alpha}|^{2}}(f^{i}_{\alpha\gamma}+h_{\alpha\overline{i},\gamma}\lambda_{\alpha}\lambda_{\gamma})
(\overline{f^{i}_{\alpha\delta}}+h_{i\overline{\alpha},\overline{\delta}}\overline{\lambda_{\alpha}}\overline{\lambda_{\delta}})\nonumber\\
&\,\,\,\,\,\,-\sum_{\alpha=1}^{\ell}R^{N}_{\gamma\overline{\delta}\alpha\overline{\alpha}}\lambda_{\gamma}\overline{\lambda_{\delta}}\nonumber
\end{align}
Putting all the above together we can get the formula (1.1).
\eproof

\bproof [Proof of Theorem \ref{tk2}]
Here we will employ the method of lemma 4.1 in \cite{YZ} which is a slightly simpler proof. Let $V_{\ell}=\sum_{\alpha=1}^{\ell}\sum_{i=1}^{n}f^{i}_{\alpha}dz^{\alpha}\otimes
e_{i}\in\Gamma(M,E)$, $E=T^{\prime\ast}M\otimes f^{\ast}(T^{\prime}N)$, $e_{i}=f^{\ast}\frac{\partial}{\partial\omega^{i}}$. Since $f$ is holomorphic map, $V_{\ell}$ is
a holomorphic section of $E$. Cleraly, $U_{\ell}=|V_{\ell}|^{2}=\sum_{\alpha,\beta=1}^{\ell}g^{\alpha\overline{\beta}}A_{\alpha\overline{\beta}}$. Thus by Bochner's formula,
we have
\begin{align}
\partial\overline{\partial}|V_{\ell}|^{2}=\langle\nabla^{\prime}V_{\ell},\nabla^{\prime}V_{\ell}\rangle-\langle\Theta^{E}V_{\ell},V_{\ell}\rangle \nonumber
\end{align}
where $\Theta^{E}$ is the curvature of the vector bundle $E$ with respect to the induced metric. Since
\begin{align}
\Theta^{E}=\Theta^{T^{\prime\ast}M}\otimes Id_{f^{\ast}(T^{\prime}N)}+Id_{T^{\prime\ast}M}\otimes f^{\ast}(\Theta^{T^{\prime}N})\nonumber
\end{align}
More precisely, we can assume that
\begin{align}
\Theta^{T^{\prime\ast}M}=-\sum_{\alpha,\beta,\delta,\eta=1}^{m}(R^{M})^{\delta}_{\alpha\overline{\beta}\eta}dz^{\alpha}\wedge d\overline{z}^{\beta}
\otimes\frac{\partial}{\partial z_{\delta}}\otimes dz_{\eta}\nonumber\\
f^{\ast}(\Theta^{T^{\prime}N})=\sum_{i,j,k,l=1}^{n}\sum_{\alpha,\beta=1}^{m}(R^{N})_{i\overline{j}k}^{l}f^{i}_{\alpha}\overline{f^{j}_{\beta}}
dz^{\alpha}\wedge d\overline{z}^{\beta}\otimes e^{\ast}_{k}\otimes e_{l}. \nonumber
\end{align}
where $(R^{M})_{\alpha\overline{\beta}\eta}^{\delta}=\sum_{\xi=1}^{m}g^{\delta\overline{\xi}}R^{M}_{\alpha\overline{\beta}\eta\overline{\xi}}$
and $(R^{N})_{i\overline{j}k}^{l}=\sum_{s=1}^{n}h^{l\overline{s}}R^{N}_{i\overline{j}k\overline{s}}$.

Hence, at the point $p$ assumed in the theorem condition, we have
\begin{align}
\langle\Theta^{E}V_{\ell},V_{\ell}\rangle=(-\sum_{\delta=1}^{\ell}R^{M}(\alpha,\overline{\beta},\delta,\overline{\delta})|\lambda_{\delta}|^{2}
+\sum_{\gamma=1}^{\ell}R^{N}(\alpha,\overline{\beta},\gamma,\overline{\gamma})|\lambda_{\gamma}|^{2}\lambda_{\alpha}\overline{\lambda_{\beta}})dz^{\alpha}\wedge
d\overline{z}^{\beta} \nonumber
\end{align}
Putting all the above together we can get the formula (1.2).
\eproof

\section{Applications}
Since in general $\|\wedge^{\ell}\partial f\|_{0}^{2}$ and $\sigma_{\ell}$ are not smooth, we will consider that
$W_{\ell}(x)$ serves a smooth barrier for $\|\wedge^{\ell}\partial f\|_{0}(x)$ and $U_{\ell}$ serves a smooth barrier
for $\sigma_{\ell}$. As stated in the previous section, The Hermitian form
$f^{\ast}h=A_{\alpha\overline{\beta}}dz^{\alpha}\wedge d\overline{z}^{\beta}$ with $A_{\alpha\overline{\beta}}=\sum_{i,j=1}^{n}f^{i}_{\alpha}\overline{f^{j}_{\beta}}h_{i\overline{j}}$. we assume that $g_{\alpha\overline{\beta}}=\delta_{\alpha\overline{\beta}}$ at $p\in M$, $h_{i\overline{j}}=\delta_{i\overline{j}}$ at $f(p)$.
After a constant unitary change of coordinates $z$ and $\omega$, at $p$, we have $df(\frac{\partial}{\partial z^{\alpha}})=\sum_{i=1}^{n}\lambda_{\alpha}\delta_{i\alpha}\frac{\partial}{\partial \omega^{i}}$ with
$|\lambda_{1}|\geq\cdot\cdot\cdot\geq|\lambda_{\alpha}|\geq\cdot\cdot\cdot|\lambda_{m}|$. For the local Hermitian metric $A=(A_{\alpha\overline{\beta}})$
and $G=(g_{\alpha\overline{\beta}})$ we denote $A_{\ell}$ and $G_{\ell}$ be the upper-left $\ell\times\ell$ blocks of them. It is easy to see
that $|\lambda_{1}|^{2}\geq\cdot\cdot\cdot\geq |\lambda_{m}|^{2}$ are the eigenvalues of $A$ (with respect to $g$). we start with some needed algebraic results (see \cite{Ni2}).
\bproposition [\cite{Ni2}, Proposition 2.1]\label{tk6}
For any $1\leq\ell\leq m$ the following holds:
\begin{align}
\sigma_{\ell}=\sum_{\alpha=1}^{\ell}|\lambda_{\alpha}|^{2}\geq\sum_{\alpha,\beta=1}^{\ell}g^{\alpha\overline{\beta}}A_{\alpha\overline{\beta}}=U_{\ell}
\end{align}
\eproposition
\bproposition [\cite{Ni2}, Proposition 2.2]\label{tk6}
For any $1\leq\ell\leq m$ the following holds:
\begin{align}
\|\wedge^{\ell}\partial f\|_{0}^{2}=\Pi_{\alpha=1}^{\ell}|\lambda_{\alpha}|^{2}\geq\frac{\det(A_{\ell})}{\det(G_{\ell})}=W_{\ell}
\end{align}
\eproposition
\bproof [Proof of Theorem \ref{tk3}]
To prove part (a), let $D=\frac{(f^{\ast}h)^{m}}{g^{m}}=W_{m}$. Since $M$ is compact, $D$ attains its
maximum at some point $p$. We assume that at $p$, $D$ is not equal to zero. Then in a neighborhood of $p$, $D\neq 0$.
The maximum principle then implies that at $p$, with respect to the coordinates specified in Theorem \ref{tk1}
\begin{align}
0\geq\Delta\log D&=S^{M}-\sum_{\delta=1}^{m}Ric_{m}^{(1)}(\delta,\overline{\delta})|\lambda_{\delta}|^{2}+
\sum_{\alpha=1}^{\ell}\sum_{i=\ell+1}^{n}\frac{1}{|\lambda_{\alpha}|^{2}}|f^{i}_{\alpha\delta}+h_{\alpha\overline{i},\delta}\lambda_{\alpha}\lambda_{\delta}|^{2}\nonumber\\
&\geq -K+\kappa(\sum_{\delta=1}^{m}|\lambda_{\delta}|^{2})\nonumber\\
&\geq-K+mD^{\frac{1}{m}}\kappa\nonumber
\end{align}
where $Ric_{m}^{(1)}$ denotes the $m$-first Ricci curvature of $(N,h)$. The above inequality implies the result.

For (b), Clearly $\|\wedge^{\ell}\partial f\|^{2}_{0}$ attains a maximum somewhere at $p$ in $M$. We assume that the coordinates at $p$ and $f(p)$
satisfy condition of theorem \ref{tk1}. Since $(M,g)$ is K\"{a}hler, we can also assume that at $p$,
$g_{\alpha\overline{\beta},\delta}=g_{\alpha\overline{\beta},\overline{\delta}}=0$, $\forall$ $1\leq\alpha,\beta,\delta\leq m$. Then we have
$\|\wedge^{\ell}\partial f\|^{2}_{0}(p)=W_{\ell}(p)$, and $W_{\ell}(x)\leq\|\wedge^{\ell}\partial f\|^{2}_{0}(x)\leq\|\wedge^{\ell}\partial f\|^{2}_{0}(p)=W_{\ell}(p)$
for $x$ in the small neighborhood of $p$. Hence, $W_{\ell}$ also attains its maximum at $p$. Now at $p$, we get
\begin{align}
0&\geq \sum_{\gamma=1}^{\ell}\frac{\partial^{2}}{\partial z^{\gamma}\partial \overline{z}^{\gamma}}(\log(W_{\ell}))\nonumber\\
&\geq S_{\ell}^{M}-\sum_{\gamma=1}^{\ell}Ric_{\ell}^{(1)}(\gamma,\overline{\gamma})|\lambda_{\gamma}|^{2}\nonumber\\
&\geq-K+\ell (W_{\ell})^{\frac{1}{\ell}}\kappa \nonumber
\end{align}
where $Ric_{\ell}^{(1)}$ denotes the $\ell$-first Ricci curvature of $(N,h)$.

The proof of part (c) is similar. We assume that $\sigma_{\ell}$ attains a maximum at $p$. We also assume that the coordinates at $p$ and $f(p)$
satisfy condition of theorem \ref{tk2}, then $U_{\ell}(x)\leq\sigma_{\ell}(x)\leq\sigma_{\ell}(p)=U_{\ell}(p)$ for  $x$ in the small neighborhood of $p$.
So, at $p$, we have
\begin{align}
0&\geq\sum_{\alpha=1}^{\ell}\frac{1}{2}(\nabla_{\alpha}\nabla_{\overline{\alpha}}+\nabla_{\overline{\alpha}}\nabla_{\alpha})U_{\ell}\nonumber\\
&\geq \sum_{\delta=1}^{\ell}Ric_{\ell}^{(2)}(\delta,\overline{\delta})|\lambda_{\delta}|^{2}-\sum_{\alpha,\gamma=1}^{\ell}
R^{N}(\alpha,\overline{\alpha},\gamma,\overline{\gamma})|\lambda_{\alpha}|^{2}|\lambda_{\gamma}|^{2}\nonumber\\
&\geq-K\sum_{\delta=1}^{\ell}|\lambda_{\delta}|^{2}+\kappa\sum_{\alpha=1}^{\ell}|\lambda_{\alpha}|^{4}\nonumber\\
&\geq-K U_{\ell}(p)+\frac{\kappa}{\ell}U_{\ell}^{2}(p)\nonumber
\end{align}
where $Ric_{\ell}^{(2)}$ denotes the $\ell$-second Ricci curvature of $(M,g)$.
\eproof

Combining theorem \ref{tk1}, theorem \ref{tk2} and theorem \ref{tk3}, we can now easily prove theorem \ref{tk3.1}.
\bproof [Proof of Theorem \ref{tk3.1}]
if $f$ is not degenerate, then $D=\frac{(f^{\ast}h)^{m}}{g^{m}}=W_{m}$ has a nonzero maximum somewhere at $p$. By using the coordinates
around $p$ and $f(p)$ specified as in the above theorem $1.4$, at $p$, we have that
\begin{align}
0\geq \Delta\log D\geq S^{M}-\sum_{\delta=1}^{m}Ric_{m}^{(1)}(\delta,\overline{\delta})|\lambda_{\delta}|^{2}\nonumber
\end{align}
This leads to a contradiction under the assumptions either $S^{M}\geq 0$ and manifold $(N,h)$ has $Ric_{m}^{(1)}<0$, or $S^{M}>0$ and
$Ric_{m}^{(1)}\leq 0$. For the second part of (a), we introduce the operator:
\begin{align}
\Psi=\sum_{\gamma=1}^{m}\frac{1}{2|\lambda_{\gamma}|^{2}}(\nabla_{\gamma}\nabla_{\overline{\gamma}}+\nabla_{\overline{\gamma}}\nabla_{\gamma})\nonumber
\end{align}
Since at $p$, $D\neq 0$, the above operator is well defined in a small neighborhood of $p$. Then, at $p$,
\begin{align}
0\geq \Psi(\log D)\geq \sum_{\gamma=1}^{m}\frac{1}{|\lambda_{\gamma}|^{2}}Ric^{M}(\gamma,\overline{\gamma})-S_{m}^{N}\nonumber
\end{align}
The above also induces a contradiction under either $Ric^{M}\geq0$ and $S_{m}^{N}<0$, or  $Ric^{M}>0$ and $S_{m}^{N}\leq0$.

The proof of (b) is similar to (a). It is worth noting that in the second part of (b), we need to introduce the following operator:
\begin{align}
\Psi_{\ell}=\sum_{\gamma=1}^{\ell}\frac{1}{2|\lambda_{\gamma}|^{2}}(\nabla_{\gamma}\nabla_{\overline{\gamma}}+\nabla_{\overline{\gamma}}\nabla_{\gamma})\nonumber
\end{align}

For (c), if $f$ is not constant, $\sigma_{\ell}$ will attains a maximum somewhere, say at $p$ and $\sigma_{\ell}(p)>0$. By using the coordinates
around $p$ and $f(p)$ specified as in the above theorem $1.4$, at $p$, we have that
\begin{align}
0&\geq\sum_{\alpha=1}^{\ell}\frac{1}{2}(\nabla_{\alpha}\nabla_{\overline{\alpha}}+\nabla_{\overline{\alpha}}\nabla_{\alpha})U_{\ell}\nonumber\\
&\geq \sum_{\delta=1}^{\ell}Ric_{\ell}^{(2)}(\delta,\overline{\delta})|\lambda_{\delta}|^{2}-\sum_{\alpha,\gamma=1}^{\ell}
R^{N}(\alpha,\overline{\alpha},\gamma,\overline{\gamma})|\lambda_{\alpha}|^{2}|\lambda_{\gamma}|^{2}\nonumber
\end{align}
if $Ric_{\ell}^{(2)}>0$, the first term is positive, the second one is nonnegative since $\widetilde{B}^{N}\leq 0$. Hence a contradiction. The same holds
if $Ric_{\ell}^{(2)}\geq 0$ and $\widetilde{B}^{N}<0$.
\eproof

\section{an integral inequality for non-degenerate holomorphic maps}
In this section we prove theorem \ref{tk5}. The proof does not use any maximum principle argument, since the curvatures of both and target spaces may not be signed in pointwise
sense. This method is essentially derived by Y. Zhang in \cite{Zhang}.
\bproof [Proof of Theorem \ref{tk5}]
Because $f$ is non-degenerate holomorphic map, we assume that $p\in M$ with $D(p)=\frac{(f^{\ast}h)^{m}}{g^{m}}(p)>0$. We also assume that the coordinates at $p$ and $f(p)$
satisfy condition of theorem \ref{tk1}. Let $\epsilon$ an arbitrary positive constant. By the formula (1.1), at $p$, we have
\begin{align}\label{5.1}
\Delta\log(D+\epsilon)&=\frac{\Delta D}{D+\epsilon}-\frac{|\partial D|^{2}}{(D+\epsilon)^{2}} \\
&=\frac{D}{D+\epsilon}(\frac{\Delta D}{D}-\frac{|\partial D|^{2}}{D^{2}})+\frac{\epsilon|\partial D|^{2}}{D(D+\epsilon)^{2}}\nonumber\\
&=\frac{D}{D+\epsilon}\Delta\log D+\frac{\epsilon|\partial D|^{2}}{D(D+\epsilon)^{2}}\nonumber\\
&\geq \frac{D}{D+\epsilon}\{S^{M}-tr_{g}(f^{\ast}(Ric_{m}^{(1)}(h)))\} \nonumber
\end{align}
Here we note that the above inequality is independent of the choice of coordinates.

Set $V=\{x\in M \mid D=0\,\, at\,\, x\}$, which is a proper subvariety (may be empty) of $M$. Then, the above inequality (\ref{5.1})
holds on $M\setminus V$ and by continuity we know it holds on the whole $M$.

Next, we fix a $\psi\in C^{\infty}(M,\mathbb{R})$ such that $e^{\psi}g^{m}$ is a Gauduchon metric on $M$. Intergrating the above inequality with
respect to $e^{(m-1)\psi}g^{m}$ over $M$ gives
\begin{align}\label{5.1}
\int_{M}\frac{D}{D+\epsilon}\{S^{M}-tr_{g}(f^{\ast}(Ric_{m}^{(1)}(h)))\}e^{(m-1)\psi}g^{m}\leq \int_{M}\Delta\log(D+\epsilon)e^{(m-1)\psi}g^{m}=0\nonumber
\end{align}
Here we have used that
\begin{align}
\int_{M}\Delta\log(D+\epsilon)e^{(m-1)\psi}g^{m}=\int_{M}(\Delta_{e^{\psi}g}\log(D+\epsilon))(e^{\psi}g)^{m}=0\nonumber
\end{align}
Now, we can easily use the same arguments in Theorem 1.1 in \cite{Zhang} to complete the proof.

\eproof

\end{document}